\documentclass[11pt]{article} 
\usepackage{hfstyle}
\usepackage{graphicx,bbm}
\usepackage[right]{eurosym}
\usepackage[T1]{fontenc}
\usepackage{amsmath}

\usepackage[polutonikogreek,english]{babel}
\usepackage[utf8x]{inputenx}

\usepackage{hfstyle}
 
\begin{document}
\title{The story of geometry told by coins}
\author{Henryk Fuk\'s, \#22079 
      \oneaddress{
         Department of Mathematics and Statistics,
 Brock University,\\
     St. Catharines, Ontario, Canada  \\
         \email{hfuks@brocku.ca}
       }
   }

%
\Abstract{
In three articles \cite{paper47,paper58a,paper58b}  published in  CNJ in 2012 and 2016 (updated author's versions available on arxiv.org),  we discussed
some links between mathematical sciences, coin minting and numismatics. This article is  a continuation (long overdue) of this cycle. It tells the  story of 
selected important developments in the history of geometry using modern commemorative coins as a background and  illustration. 
}
\maketitle
Geometry is one of the oldest branches of mathematics. It studies properties of space
and properties of objects (figures) in space, using concepts such as distance, shape and size. 
Almost all ancient civilizations developed some sort of geometrical concepts,
especially those  concerning the measurement of lengths,  areas and volumes. The
Mesopotamian and Egyptian civilizations are believed to be the first
cultures where such concepts appeared, but the ancient Greek
mathematicians are credited to be the true originators of rigorous geometric procedures and were the first to use deductive reasoning to solve geometric problems. 
In the popular mind, Greek geometry is often associated with Pythagoras and
his famous theorem about right triangles. This was discussed in \cite{paper58b},
thus we will not repeat it here, only adding that a new coin commemorating
Pythagoras appeared since the publication of \cite{paper58b}, namely the 10 Euro
silver  coin  issued in 2013 by the Bank of Greece under the program   ``Greek culture -- mathematicians''.
\begin{center}
\includegraphics[width=6cm]{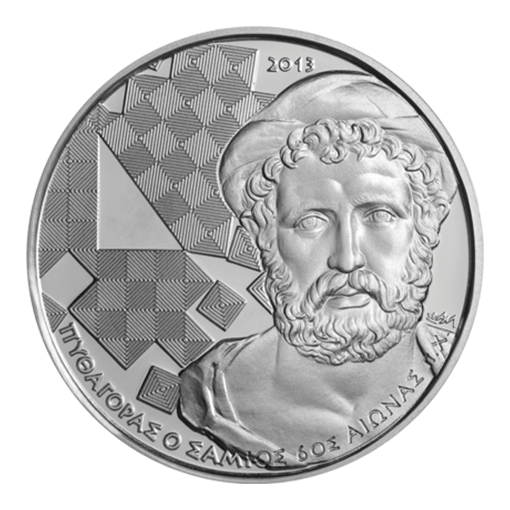}
\includegraphics[width=6cm]{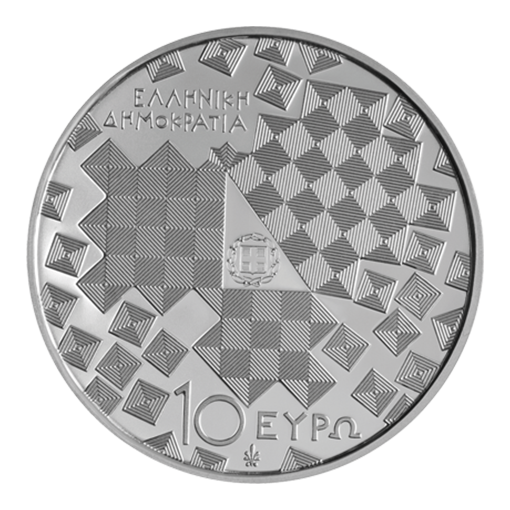}
\end{center}

 In spite of Pythagoras' popularity  among students of elementary geometry, 
 historians of mathematics tend to emphasize the achievements and influence of other ancient Greek mathematicians, often referred to as the ``big three'', namely
Euclid (325--270 BC),
Archimedes of Syracuse (c. 287--c. 212 BC),
  and  Apollonius of Perga (c. 240 BC--c. 190 BC). 
  Euclid is considered to be the father of geometry,  famous for his book \emph{Elements},  the most successful textbook of all time.  
 In terms of the number of editions, it is second only to the Bible,
and has been used as a textbook  until the early 20-th century.
  In the English language alone, there are nearly 20 different translations. 
\begin{center}
\includegraphics[width=6cm]{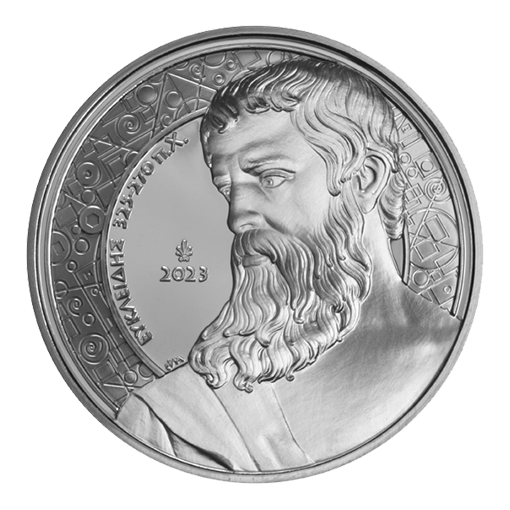}
\includegraphics[width=6cm]{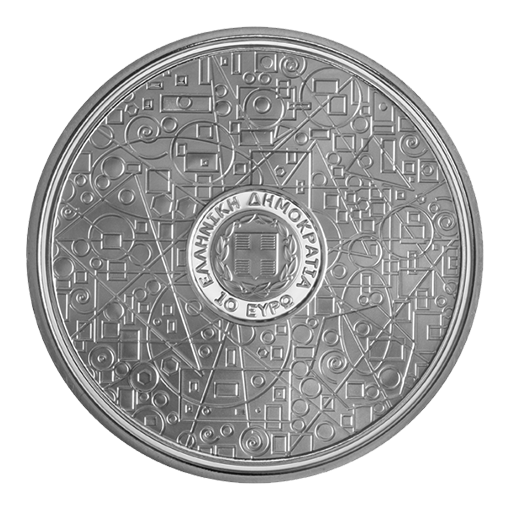}
\end{center}
 It is, therefore, somewhat surprising that Euclid had to wait nearly 23 centuries to be commemorated on a coin. This happened in 2023, when the National Bank of Greece
 issued a silver 10 Euro featuring a portrait of Euclid on the obverse,
 designed by Geórgios Stamatópoulos. The portrait
 is obviously a product of the artist's  imagination since no image of Euclid
 survived to our times. The reverse is a mosaic of various geometrical figures
 whose properties are investigated in \emph{Elements},
 such as lines, circles, ellipses and polygons. An observant reader will notice that
 spirals are present too -- this is a bit of an  anachronism as there is no mention of spirals in Euclid's works.
 
 Nevertheless, the next of the ``big three'', Archimedes, studied spirals in great details.
 He, just like Euclid, had to wait till 2015 to be featured on a coin. 
 This is also a silver proof coin, part of the same Bank of Greece series as Euclid's coin and designed by the same engraver. Archimedes, one of the greatest mathematicians
 in history, developed ingenious methods for calculating areas and volumes
 of geometric figures, including the circle, ellipse and sphere as well as the paraboloid and hyperboloid of revolution.  He invented the ``method of exhaustion'' for this purpose, anticipating modern infinitesimal calculus developed in
  the late 17th century. The bust of Archimedes on the aforementioned coin
  is shown in the background of a system of compound pulleys, one of his 
  numerous inventions, a simple machine also known as a ``block and tackle''.
The circles on the reverse, with one of them playing the role of zero
in the ``10 Euro'' inscription, likely pay homage to the alleged 
last words of Archimedes,  ``Noli turbare circulos meos!'' (``Do not disturb my circles!''). According to the Roman writer Valerius Maximus, Archimedes  uttered these words (slightly differently phrased) to a Roman soldier who came to
kill him during the siege of Syracuse in 212 BC, referring to a drawing of geometric figures which Archimedes had outlined in the sand shortly before the deadly encounter. 
\begin{center}
\includegraphics[width=6cm]{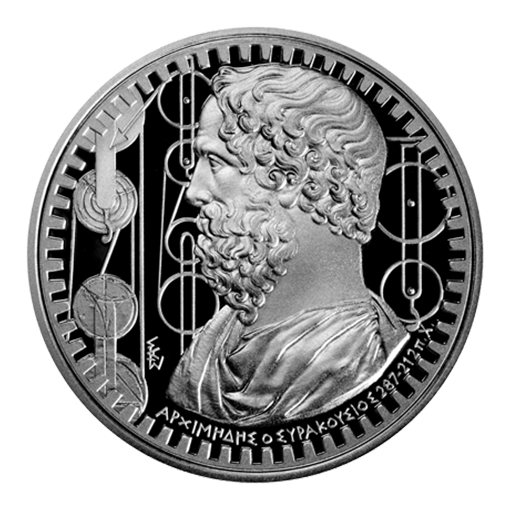}
\includegraphics[width=6cm]{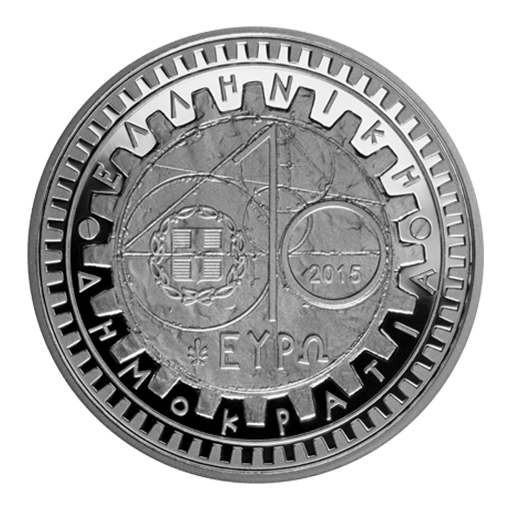}
\end{center}
 
 The third of the ``big three'', Appolonius of Perga, is known for his work on
 the so-called conic sections, curves obtained  from a cone's surface intersecting a plane. These curves are known today as Appolonius named them, namely as the ellipse, parabola and hyperbola. Appolonius does not have his own coin yet, but
 one can hope that this changes in the near future given that the  Bank of Greece continues the  program   ``Greek culture -- mathematicians''. In 2024
 another coin belonging to this program has been produced, this time
  with Thales of Miletus (624 --540 BC).
 Everyone who was fortunate enough to be exposed to  traditional geometry
instruction  in high school   certainly remembers his name thanks to the
basic proportionality theorem commonly known as the  Thales' theorem:
if BC and DE are parallel, then DE/BC = AE/AC = AD/AB. No writings of Thales
survived, but  Thales' theorem appears in the first book of Euclid's Elements as the 26th proposition.
 \begin{figure}
\begin{center}
\includegraphics[width=5cm]{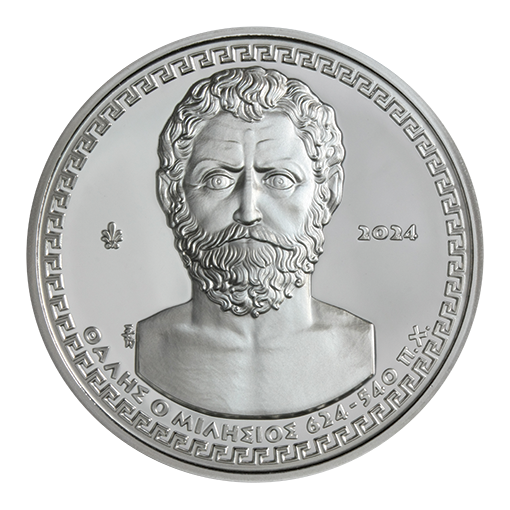}
\includegraphics[width=5cm]{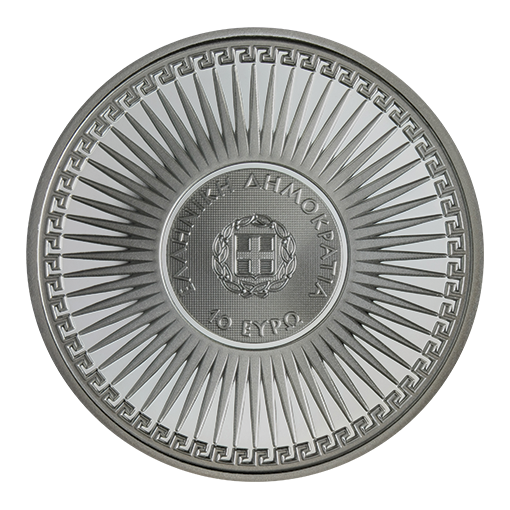}
\includegraphics[width=5cm]{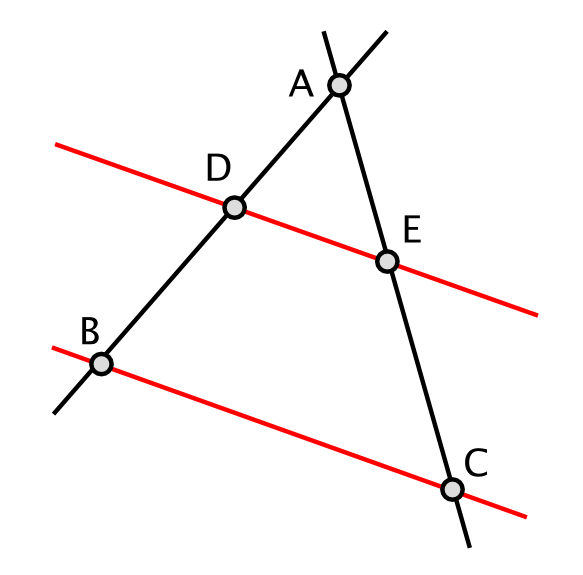}
\end{center}
 \end{figure} 
And here is my gripe: the rather enigmatic reverse of the Thales' coin is, in my opinion, a missed opportunity. The illustration of the basic proportionality theorem 
almost begs to be included.

 As we mentioned at the beginning, other ancient civilizations also 
 developed geometric ideas and made progress in expanding geometric
 knowledge. The Chinese civilization was certainly one of them.  Zu Chongzhi (429--500) was one of the most prominent
 Chinese mathematicians of late antiquity whose greatest achievement
 was calculating $\pi$ with eight digit precision. He 
 found the upper and lower bounds on the exact value of $\pi$, $3.1415926 < \pi < 3.1415927$,
 and gave the rational approximation $\pi\approx 355/113$.
\begin{figure}
\begin{center}
\includegraphics[width=6cm]{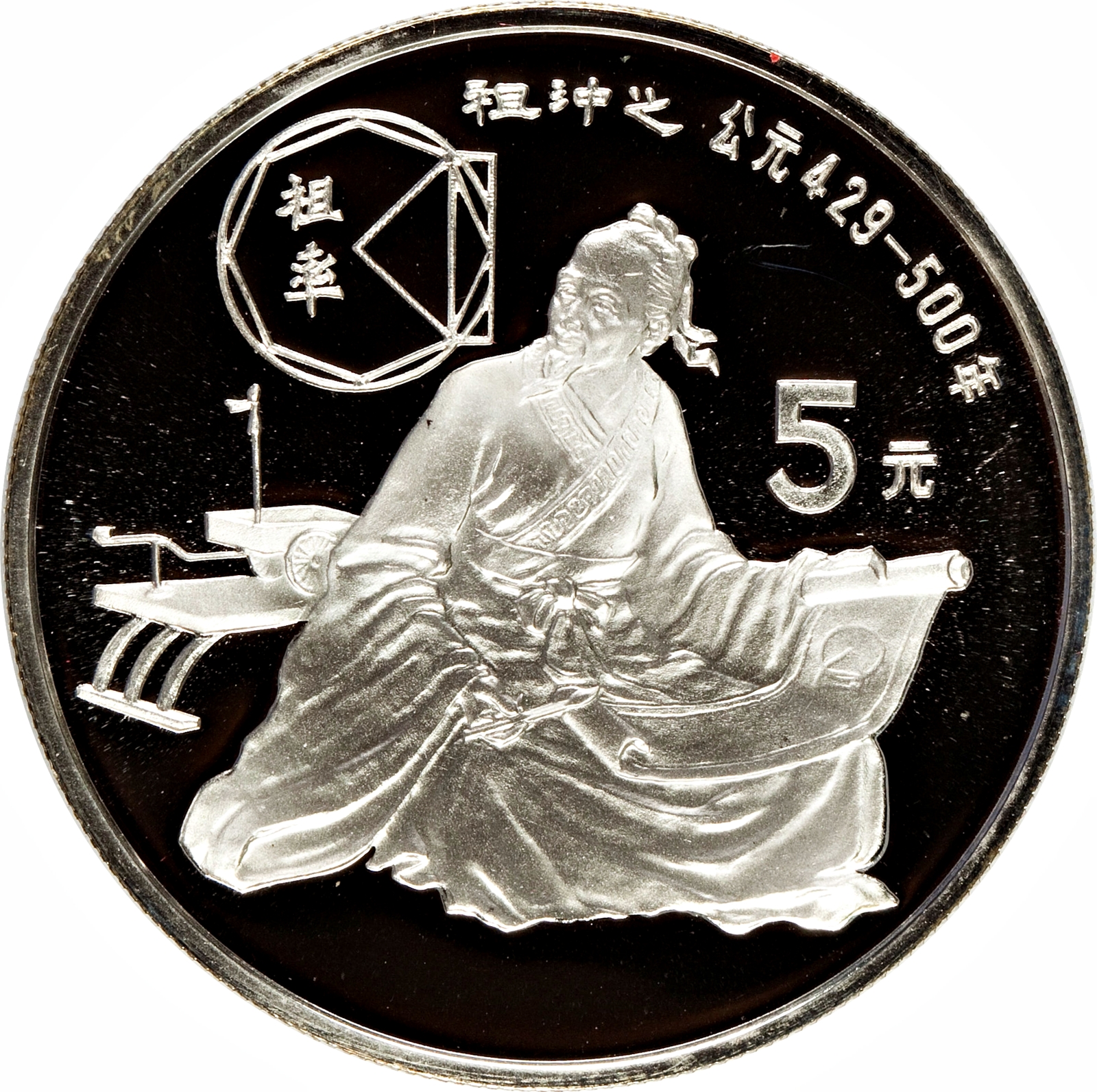}
\includegraphics[width=6cm]{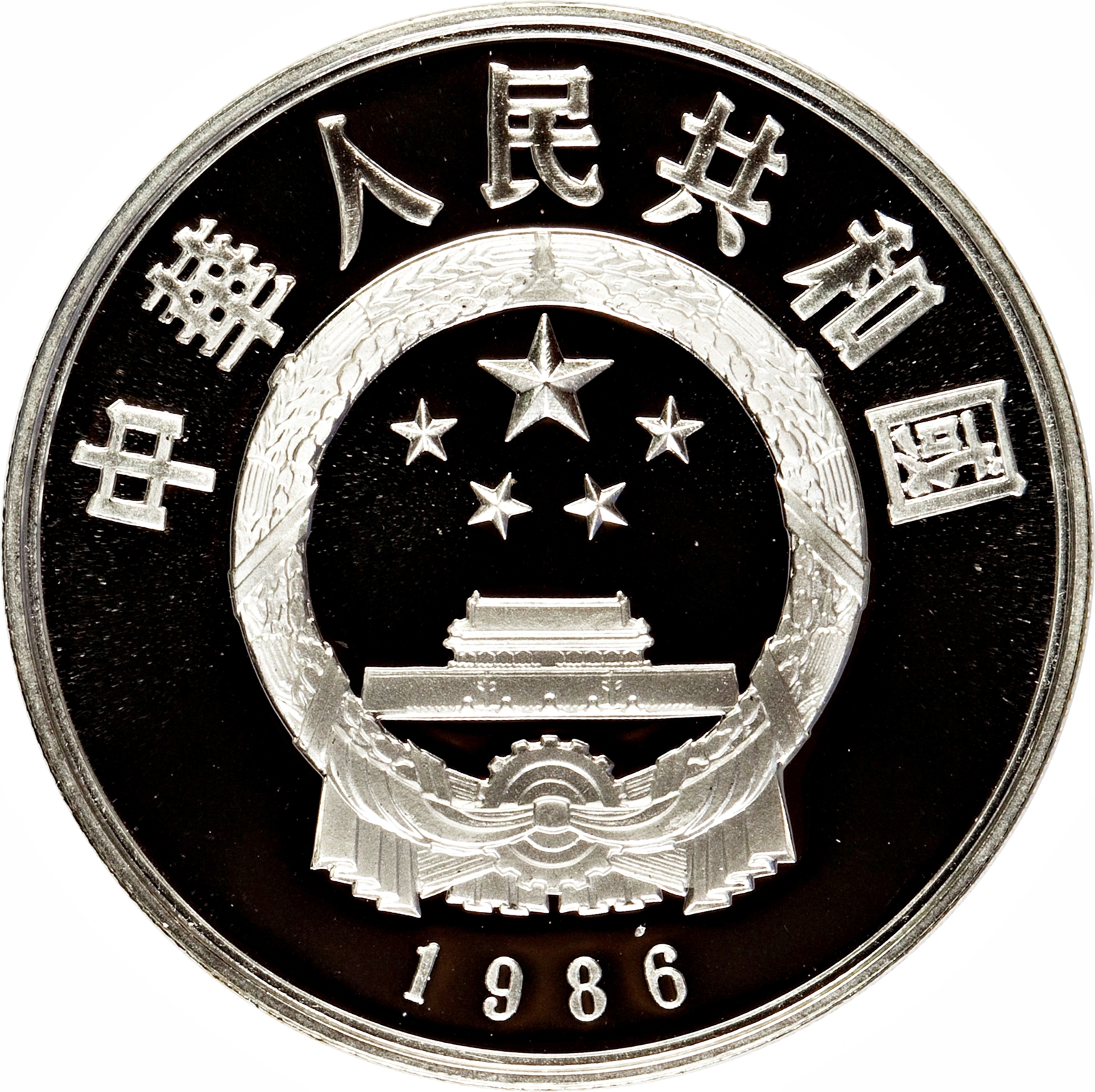}
\end{center}
 \end{figure}  
 The 5 Yuan silver coin issued in 1986 by the People's Republic of China
 presents a full-body portrait  of  Zu Chongzhi in the background of the hammer mill, one of his mechanical inventions. In the top-left part we see the polygon
 inscribed in a circle. In deriving his results regarding $\pi$,
 Zu Chongzhi approximated the circle with $n$-sided polygons, using
 increasing values of $n$. It is not known  exactly how Zu Chongzhi obtained his interval of bounds on $\pi$, but some authors conjectured that he started
 with a hexagon and then doubled the number of its sides many times, 
 obtaining  polygons with 12, 24, 48  sides, etc. Internet sources
 such as Wikipedia authoritatively claim that he 
ended up with  a polygon of  $3 \cdot 2^{12}$ sides, but there are 
valid reasons to doubt this.
 Historians of mathematics suspect that he actually started with a square,
 doubling the number of its sides until he obtained a polygon with
 $2^{15}$ sides. Fortunately, the designer of the coin did not have access to 
 Wikipedia in 1986, thus he leaned toward the more probable hypothesis,
 showing an octagon on the coin. This corresponds to the first step in the 
 approximation: 4 sides of a square, after doubling, result in a polygon with 8 sides, i.e., octagon.

Coming back to Greek mathematicians, it is often said that after the death of Archimedes  geometry suffered a decline. Indeed, Pappus of Alexandria (c. 290--c.  350) and a commentator on  Euclid's Elements Proclus Lycius (412--485) are the only Greek names in late antiquity worth mentioning as long as abstract geometry is considered.
Nevertheless, it is fair to say that practical geometry enjoyed a much better
fate. Claudius Ptolemy  (c. 100-c. 170) produced his seminal astronomical treatise now known as \emph{Almagest}, where he described his geocentric model of universe based solely on geometric principles. He also wrote a book
known as  \emph{Geographia} where he explained how to draw maps using geographical coordinates, similar to what we use today. In 1992 Cuba issued a 10 Peso silver coin
honoring Ptolemy  together with a 15-th century map maker and mathematician 
Paolo Toscanelli (1397--1482).  
\begin{center}
\includegraphics[width=6cm]{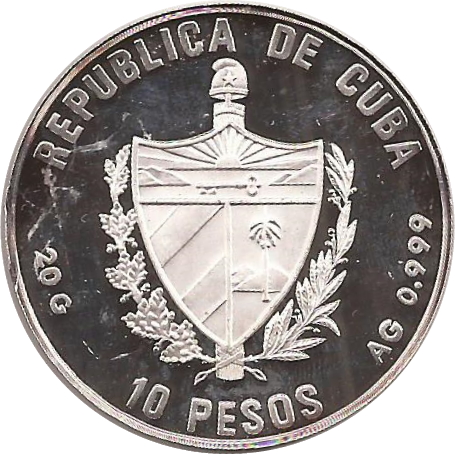}
\includegraphics[width=6cm]{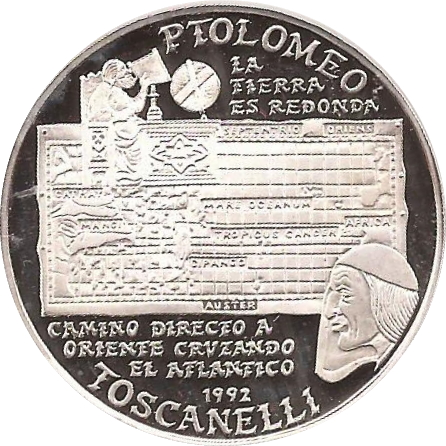}
\end{center}
The reason for including them both is
because Toscanelli used Ptolemy’s geographic information from  \emph{Geographia} 
and came up with  the idea of sailing west to reach Asia by covering a distance of only 6,500 miles.  He prepared a map illustrating his idea.
Columbus had Toscanelli’s map and it convinced
him to try to find a direct water route from Europe to Asia. A modern rendering of the map (the original is lost) is featured on the Cuban coin, which  
commemorates the 500th anniversary of Columbus' voyage. One can see from this map that
Cipango (Japan) is separated from Europe and Africa by the Mare Oceanum (Ocean Sea).
 
We will return to the subject of map making shortly, but for now let us continue the  discussion of developments in geometry in a chronological way.  In the Middle Ages
not much happened in abstract geometry, although in the late Middle Ages  practical geometry flourished again. The monumental Gothic cathedrals of Europe are a living testament of the skills of medieval architects who used various geometric methods
to design, layout and built these magnificent structures.
 
At the beginning of the 15-the century we encounter another crucial development
in  practical geometry. Filippo Brunelleschi  (1377--1446), featured on the 1000 Lira coin of San Marino from 1977 (400-th anniversary of his birth) was 
an Italian architect and designer. In the years   1415--1420 he conducted
a series of experiments which led him to the discovery of the linear perspective,
a method of representing a three dimensional object on a two-dimensional
canvas, producing the illusion of depth.
He is credited to be the first person to precisely describe  the principles of linear perspective. The San Marino coin shows on the obverse a group of people
on a raised platform or scaffolding, presumably representing 
 the architect Brunelleschi discussing something with builders. The beams perpendicular to the plane of the coin converge to a point located outside of the
 coin, following the rules of linear perspective, although with a bit of fisheye distortion added
 to harmonize with the circular shape of the coin. 
\begin{center}
\includegraphics[width=6cm]{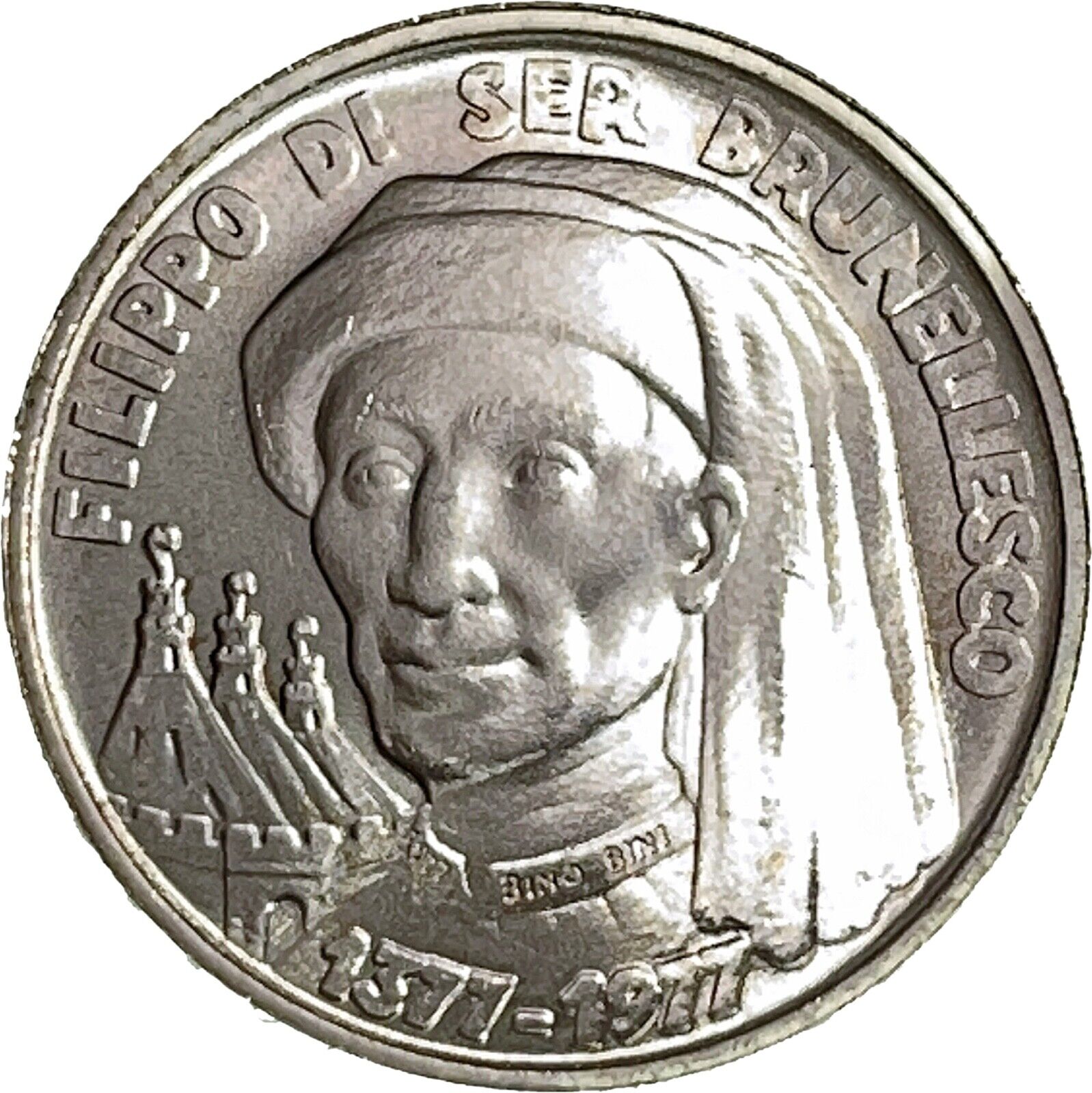}
\includegraphics[width=6cm]{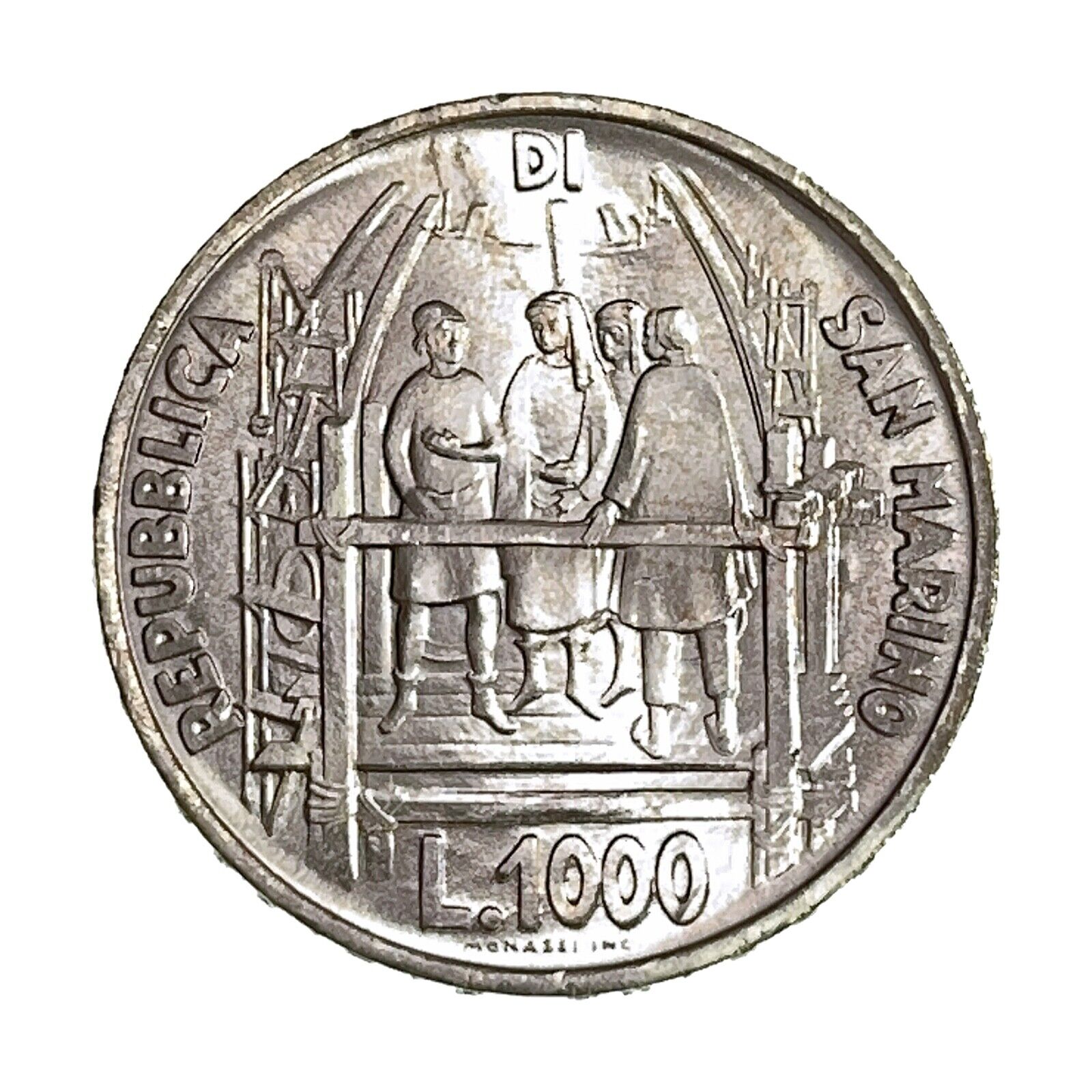}
\end{center}
In the 1470s an Italian painter  Piero della Francesca (1415--1492) further developed Brunelleschi's ideas 
and produced an influential work titled  \emph{De prospectiva pingendi} (On the Perspective of Painting). Inspired by the wok of Euclid,  Piero della Francesca
described mathematical principles of perspective and explained how to create 
perspective with colors. The Italian 500 Lire coin from 1992 shows the
portrait of della Francesca on the obverse and one of his most
famous frescos, \emph{The Resurrection}, on the reverse. While the composition of the \emph{Resurrection} is based on some geometric principles (right), the use of  perspective is better illustrated  with della Francesca's the  \emph{Flagellation of Christ} (below), where 
a number of straight lines clearly converge to a single point.
\begin{center}
\includegraphics[width=4cm]{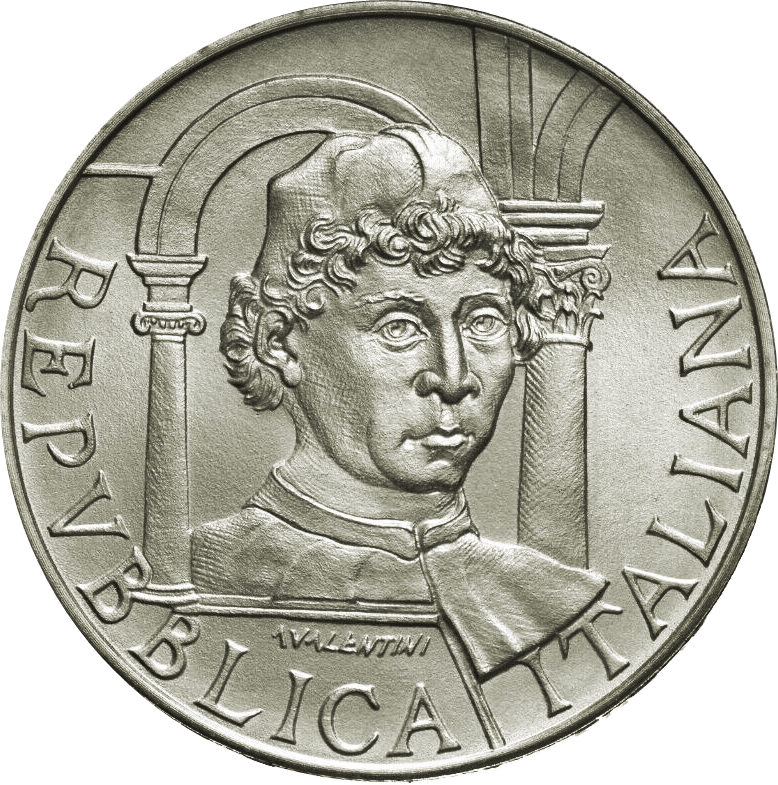}
\includegraphics[width=4cm]{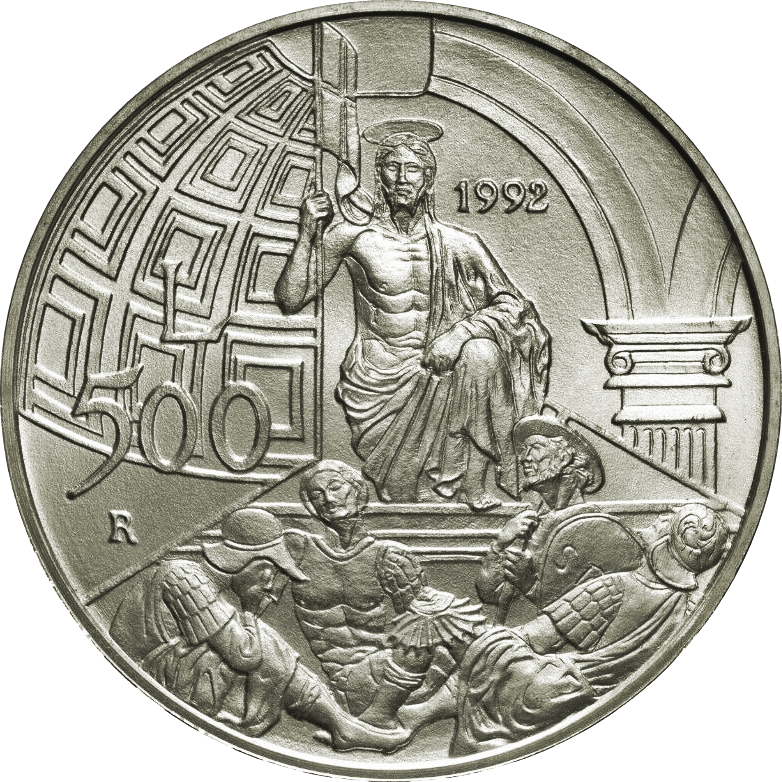}
\includegraphics[width=5cm]{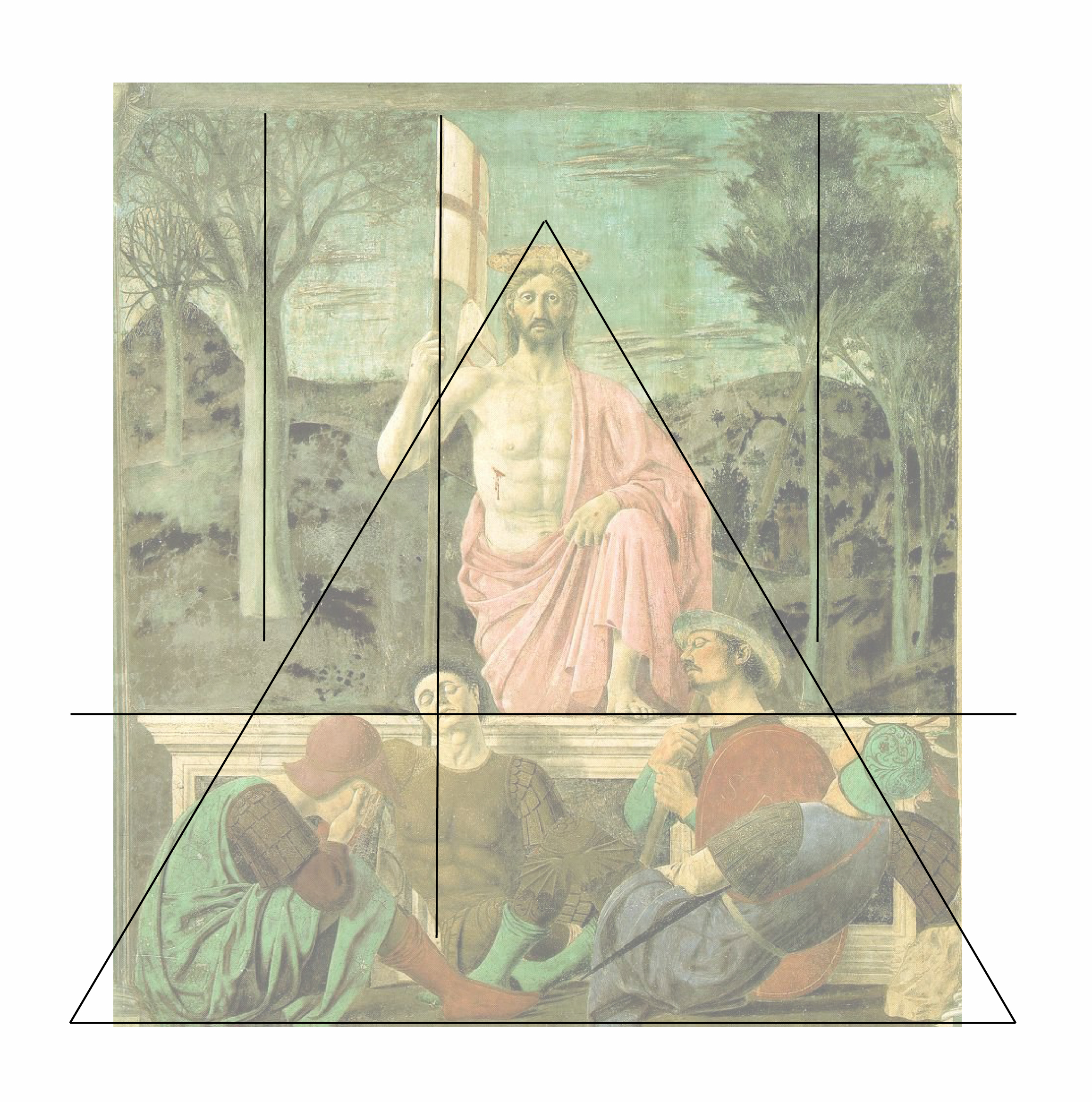}\\[1em]
\includegraphics[width=6cm]{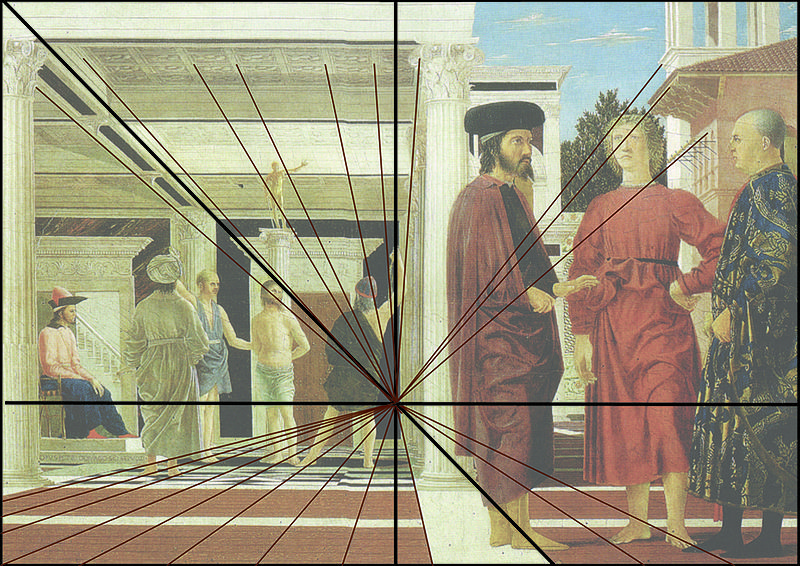}
\end{center}
After della Francesca, at the beginning of the 16th century, a Franciscan friar 
Luca Pacioli (1447--June 1517) further advanced the theory of perspective 
in his book \emph{Divina proportione} (1509). The book discusses  the notion of proportions in visual arts in the context of perspective with
applications  to architecture. The 500 Lira Italian coin from 1994 features his portrait circumscribed by the inscription ``1494 Luca Pacioli 1994''. The year
1494 represents the publication date of his other influential work \emph{Summa de arithmetica} which, among other things, contains the first published description of the double-entry bookkeeping system. This work earned him the title of ``father of accounting''.
\begin{center}
\includegraphics[width=6cm]{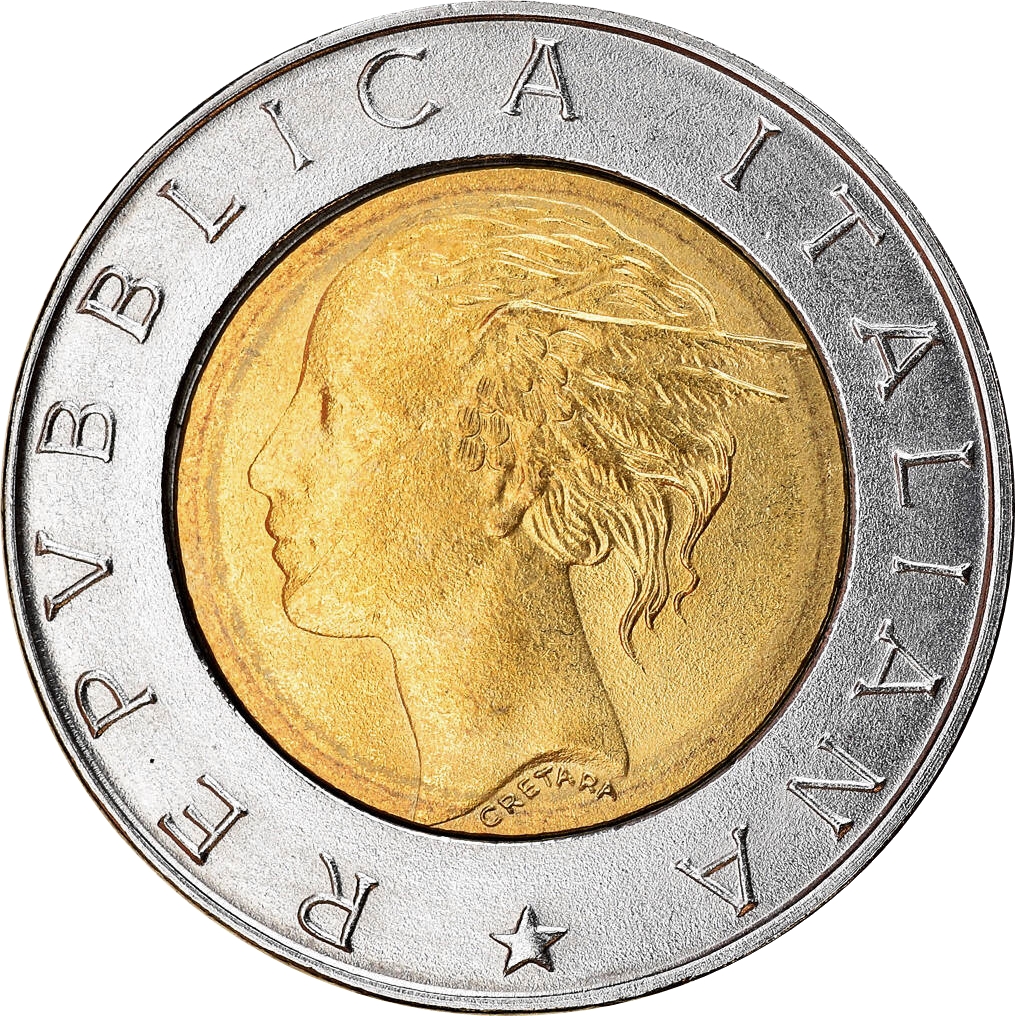}
\includegraphics[width=6cm]{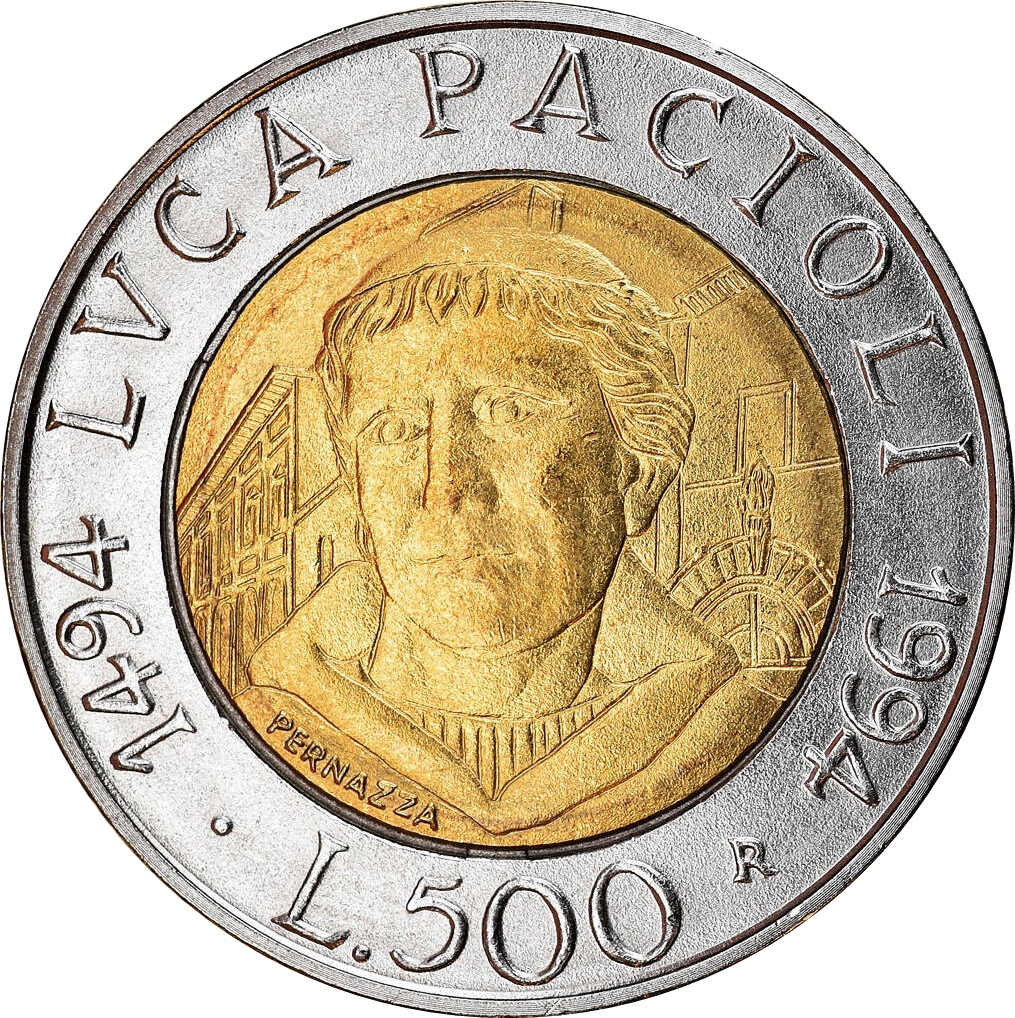}
\end{center}

After discussing the developments of practical geometry, let us return to more abstract problems. In the first half of the 17-the century an important development
in abstract geometry emerged, which could be described as the ``marriage of
geometry and algebra''. This was mainly due to René Descartes (1596--1650) who formulated the principles of analytic geometry or, in other words,  the study of geometry using a coordinate system. Today analytic geometry is often called Cartesian geometry and Descartes is credited as its sole founder, although
it  is worth stressing that Pierre de Fermat (1607--1665)
also contributed to the subject. 
The coin of France issued in 1991 (100 Francs) shows the head of Descartes
with his hair smoothly transitioning into books from the left and into
a scroll from the right. Somewhat more related to geometry is the 50 Lira coin of San Marino from 1996 with the bust of Descartes inscribed within Cartesian axes.
\begin{center}
\includegraphics[width=6cm]{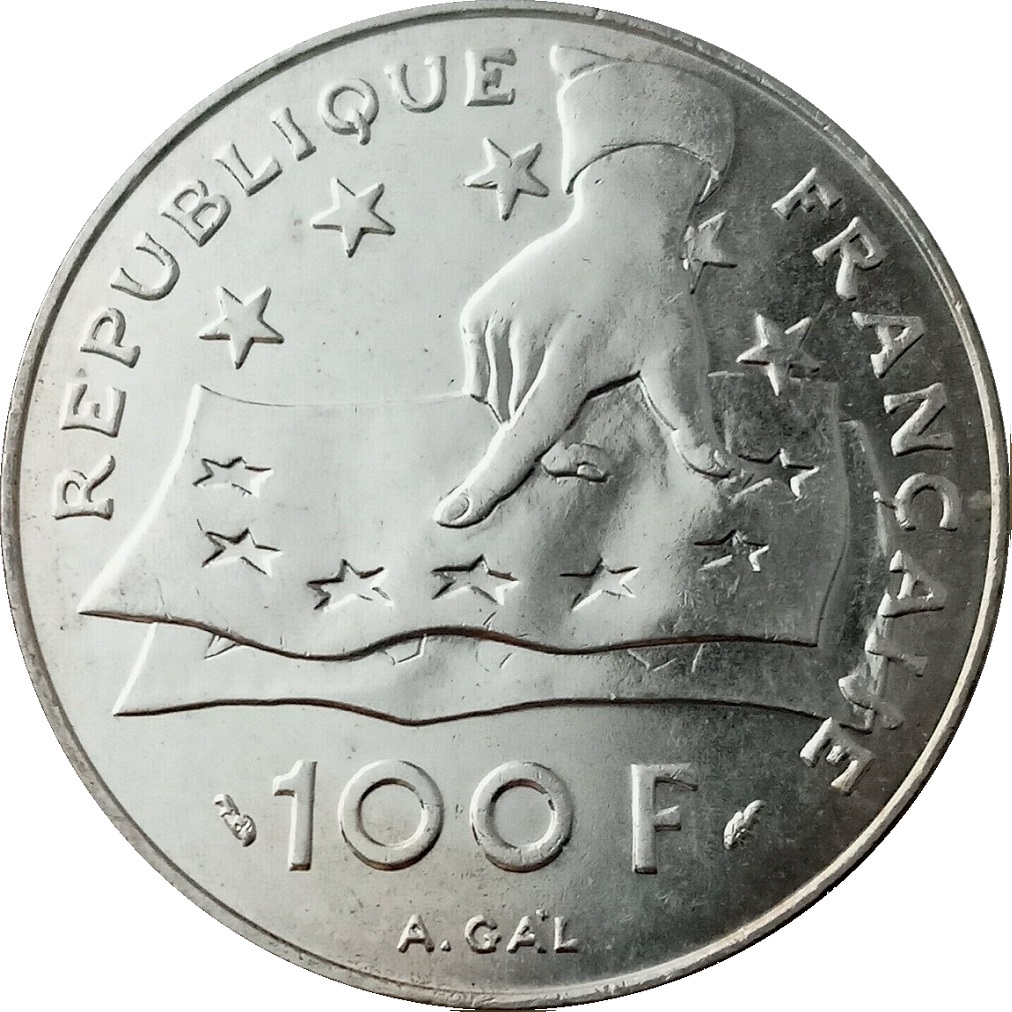}
\includegraphics[width=6cm]{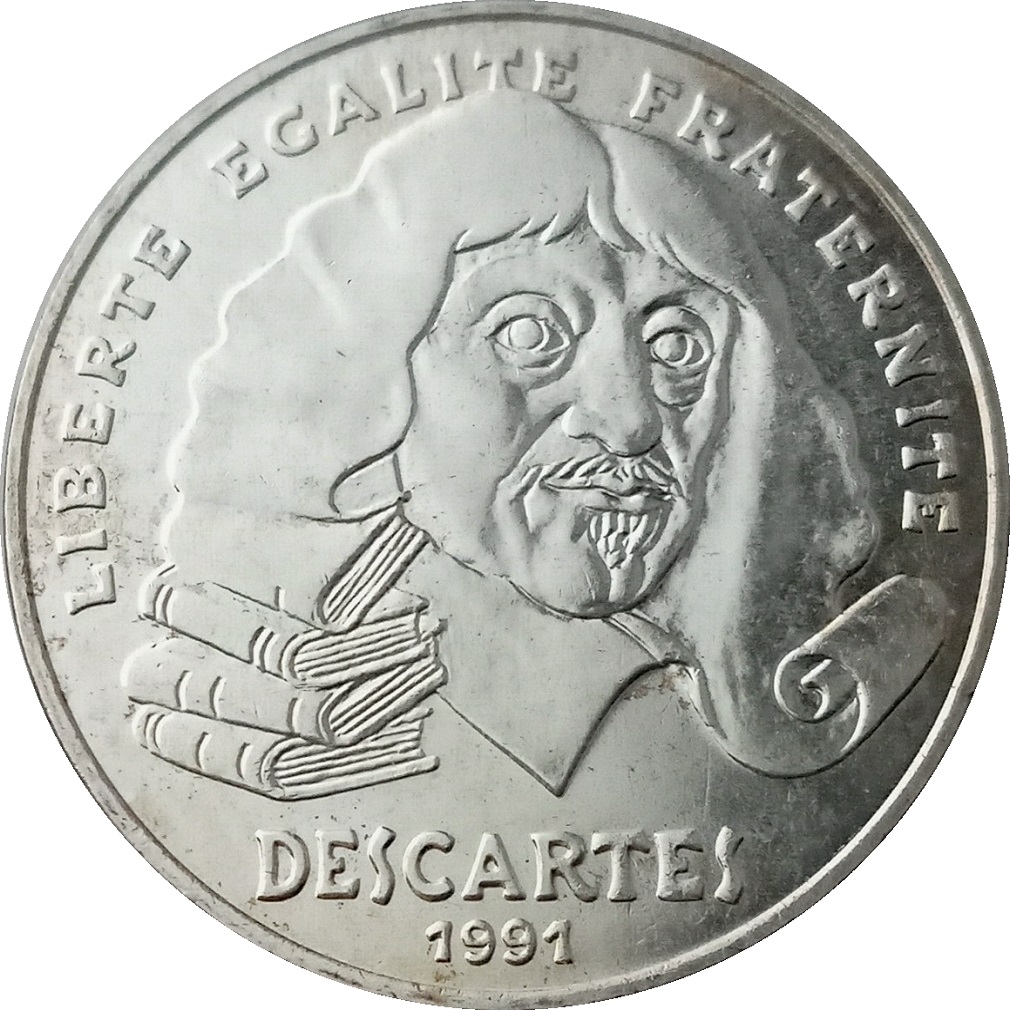}\\
\includegraphics[width=6cm]{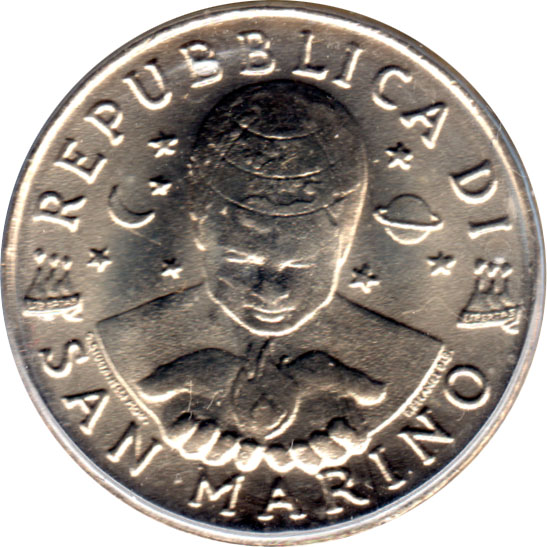}
\includegraphics[width=6cm]{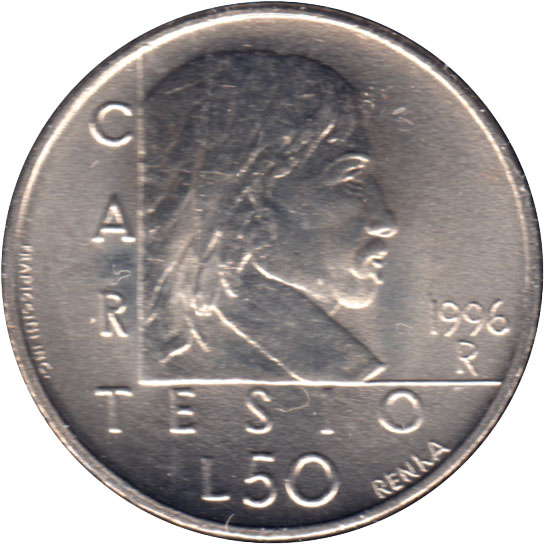}
\end{center}

In moving on to the 18th and 19th centuries we encounter the towering figure of
Carl Friedrich Gauss (1777--1855), one of the most influential mathematicians in the history of mathematics, often referred to as the ``Prince of Mathematicians''.
This is not a place to describe his countless contributions to various branches of mathematics, thus we will only mention geometry. Before we do this, let us mention that in 1977, on the 200th anniversary of Gauss' birth, two ``competing''
coins were issued by the then two distinct German states, the German Democratic Republic
and the Federal Republic of Germany. The East-German coin had the denomination
of 20 Marks (Eastern marks, of course) and carried  the
graph of the Gaussian curve, the well known  probability distribution function widely used in statistics,  also known as the ``bell curve''. This a very well chosen
symbol of Gauss' contributions to mathematics as it is instantly recognizable
and entrenched in  popular culture.
\begin{center}
\includegraphics[width=6cm]{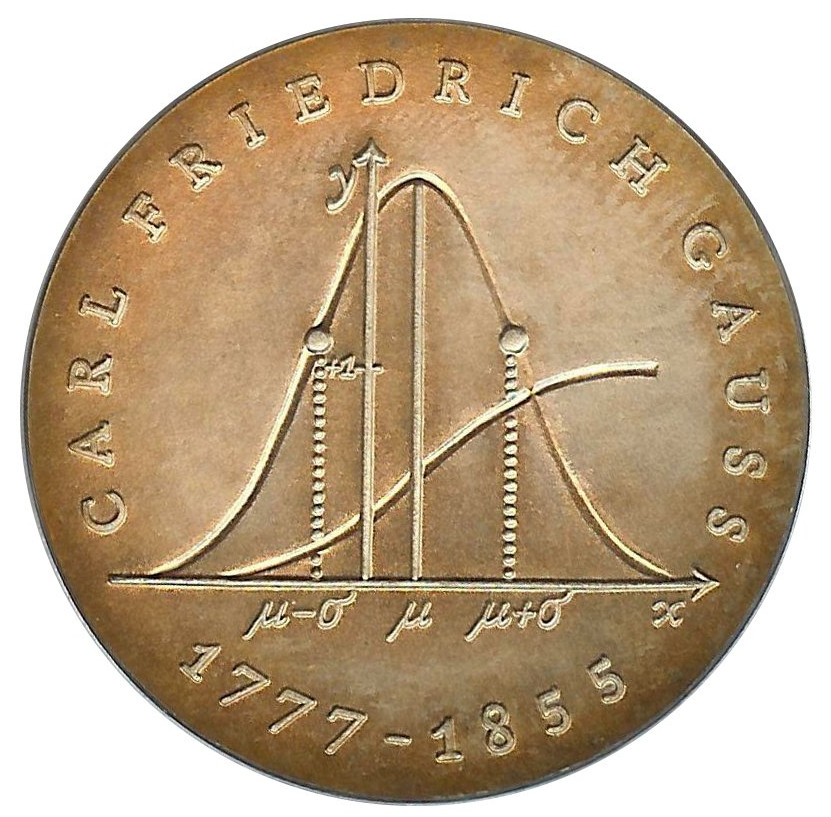}
\includegraphics[width=6cm]{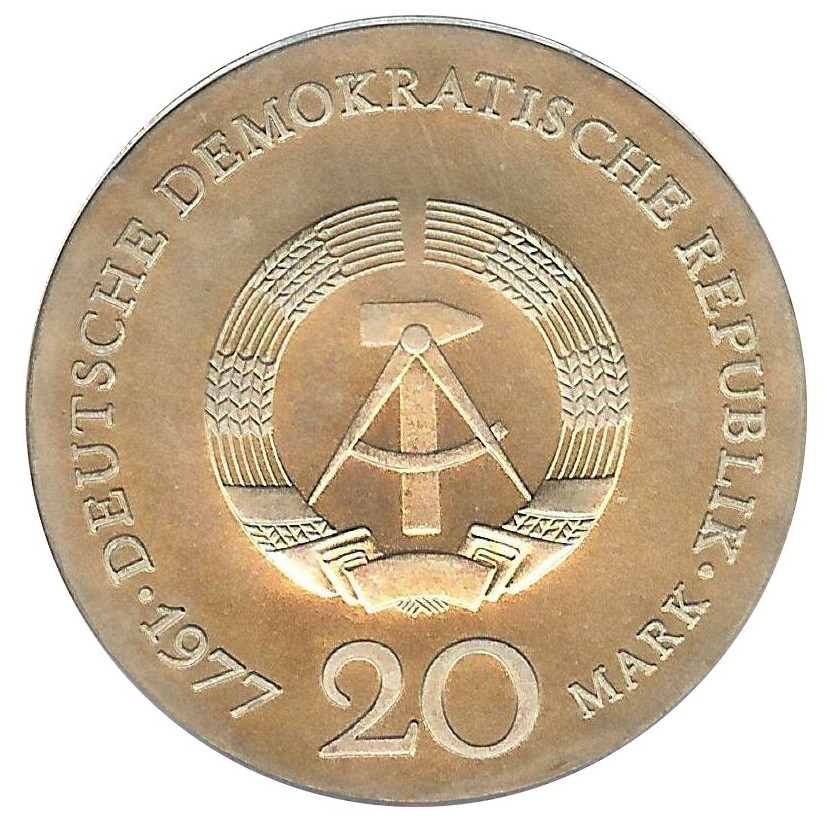}
\end{center}

The West-German coin, with a 5 DM (Deutsche Mark) denomination, is a more typical commemorative
coin with the portrait of Gauss closely inspired by the well known
painting by Christian Albrecht Jensen from 1840 (shown on the right),
with Gauss wearing his characteristic hat.
Although in 1970's the black market exchange rate was 5 to 10 eastern Marks to one Deutsche Mark, making the 5 DM coin more valuable than the 20 eastern Mark
coin of DDR, today the East German coin reaches much higher prices
on the collector market than its West German counterpart. 
Of course this mainly reflects the difference in mintage  (55,000 versus 1 million).
\begin{center}
\includegraphics[width=4cm]{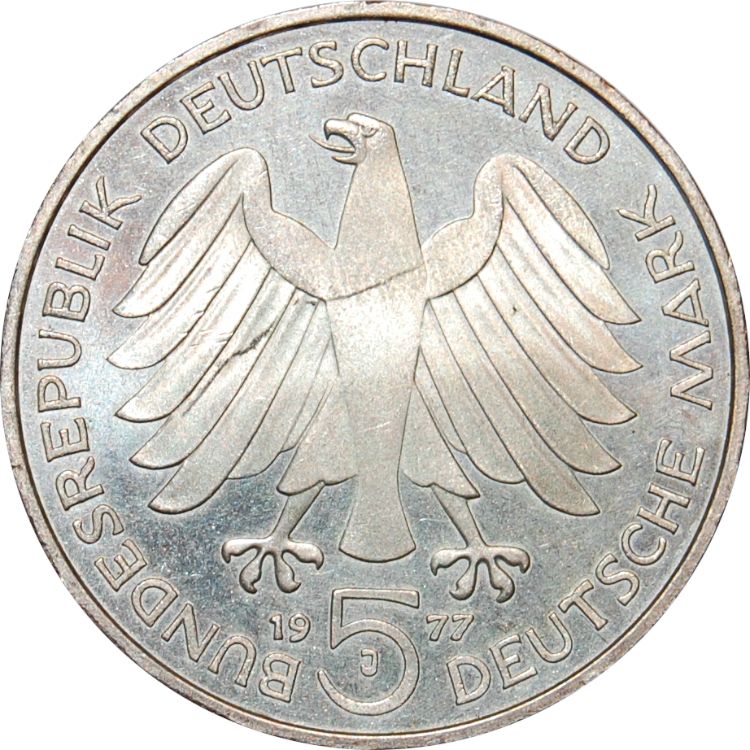}
\includegraphics[width=4cm]{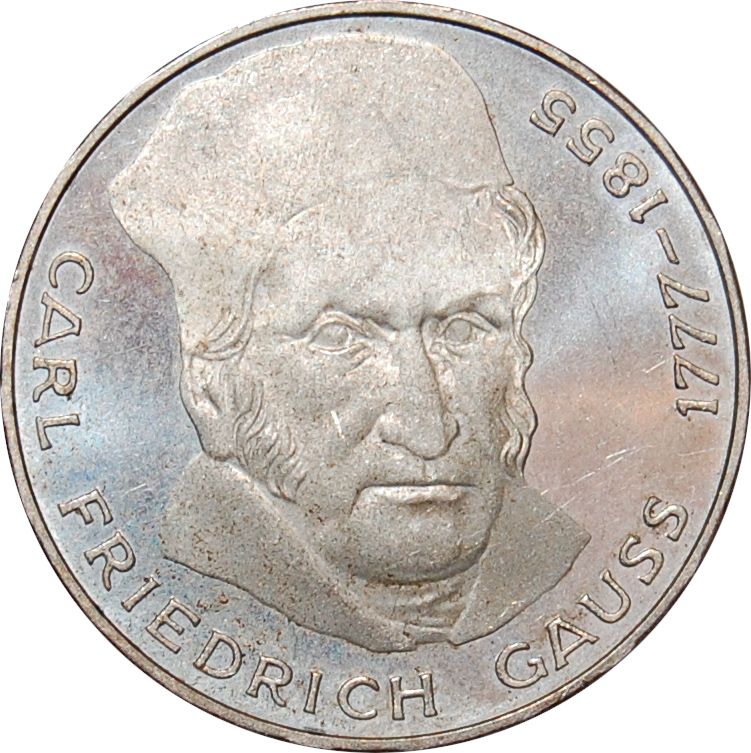}
\includegraphics[width=4cm]{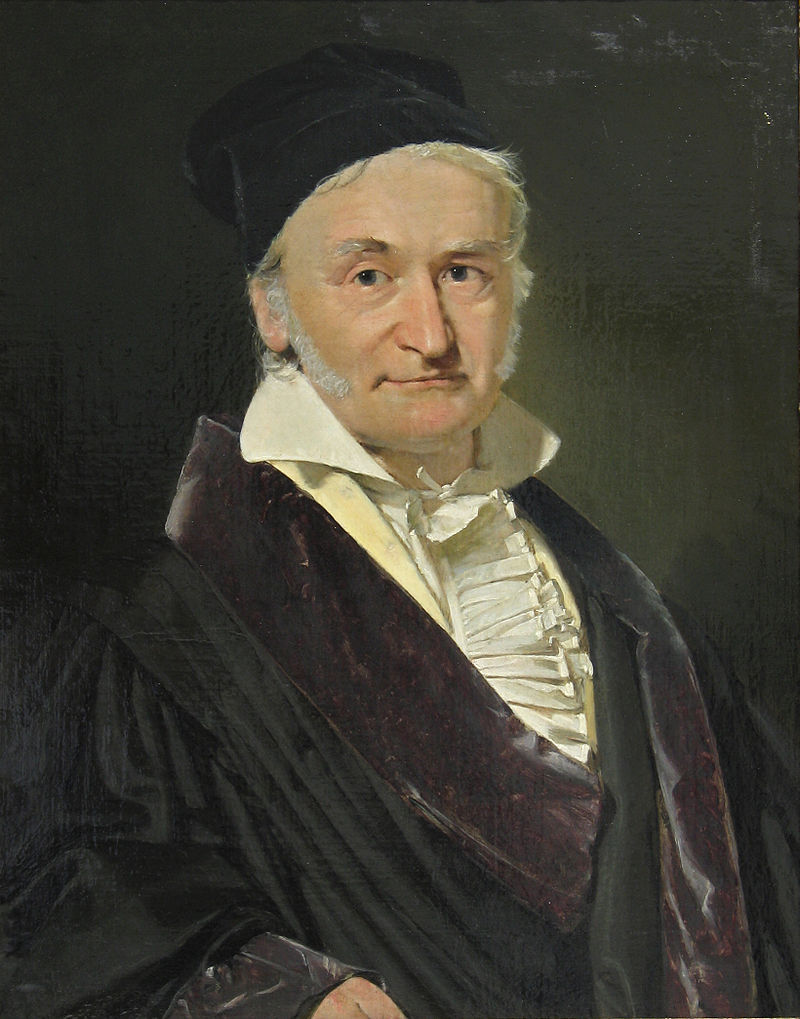}
\end{center}

Gauss is considered one of the discoverers of non-Euclidean geometry, and more about this later.
Yet Gauss also initiated differential geometry of surfaces, which is
different from Cartesian geometry because it  studies surfaces from a point of view
of an ``internal'' observer  constrained to move on the surface as if it were two-dimensional.  One of Gauss' greatest achievements in this area was his
\emph{Theorema Egregium} (remarkable theorem), which informally says that the curvature of a surface can be determined entirely by measuring angles and distances on the surface itself.

 One of the consequences of this theorem is that a 
 sphere's surface cannot be represented on a plane without distortion,
 thus one cannot produce a map of a significant area of the Earth which would be free
 of distortion. Consider, for example, a frequently used map projection known 
 as the Mercator projection, invented by a Flemish geographer Gerardus Mercator (1512 --1594), featured on the 5 Mark silver German coin of 1969.
\begin{center}
\includegraphics[width=4cm]{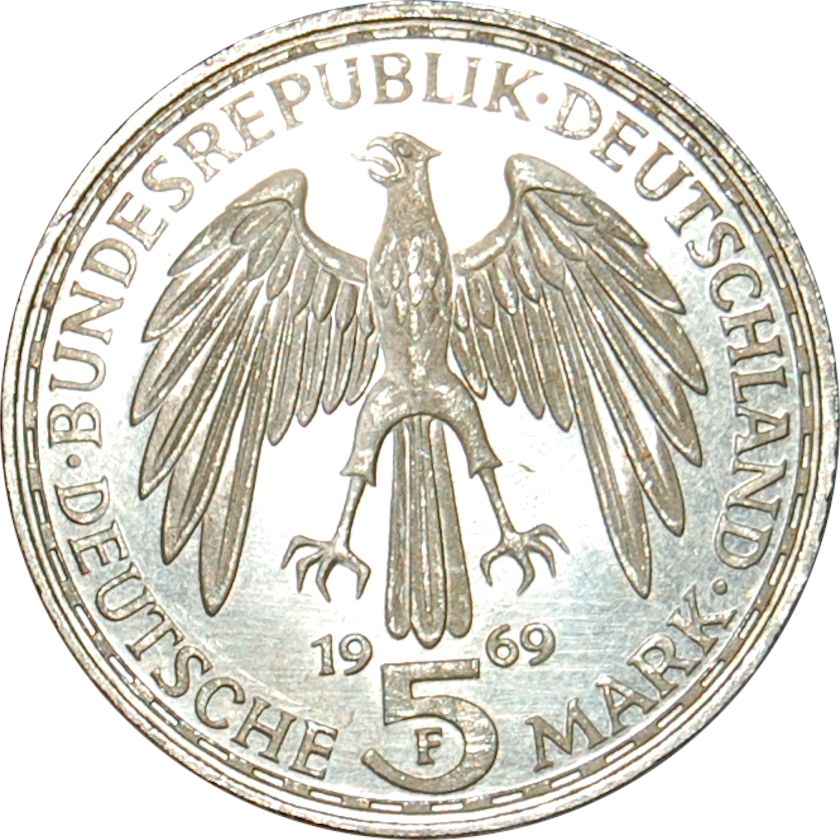}
\includegraphics[width=4cm]{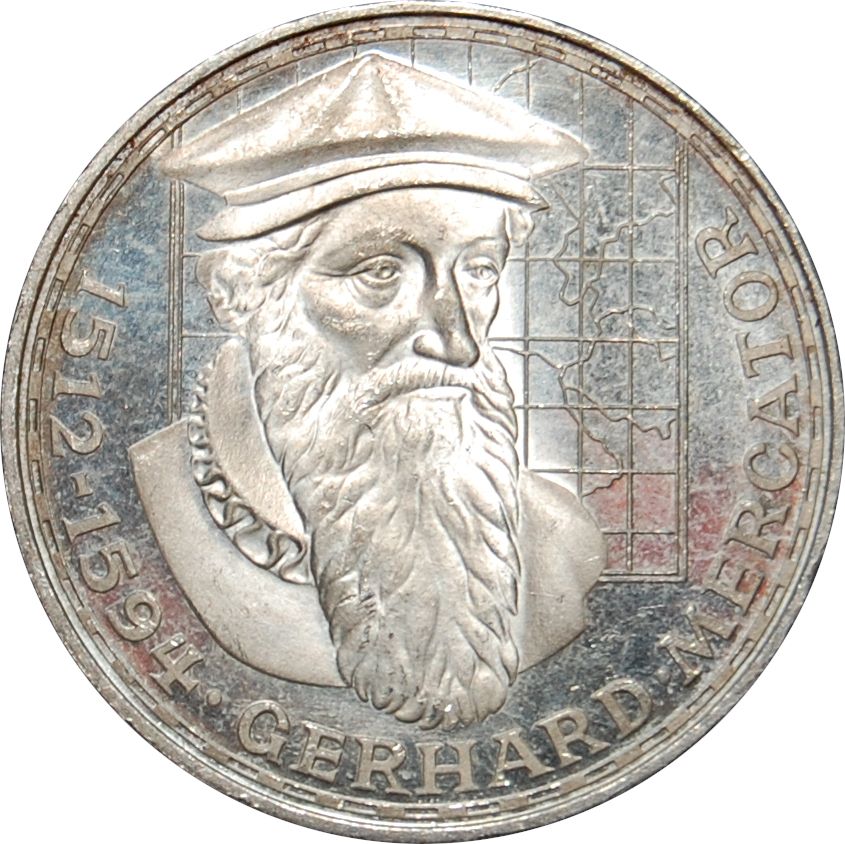}
\includegraphics[width=5cm]{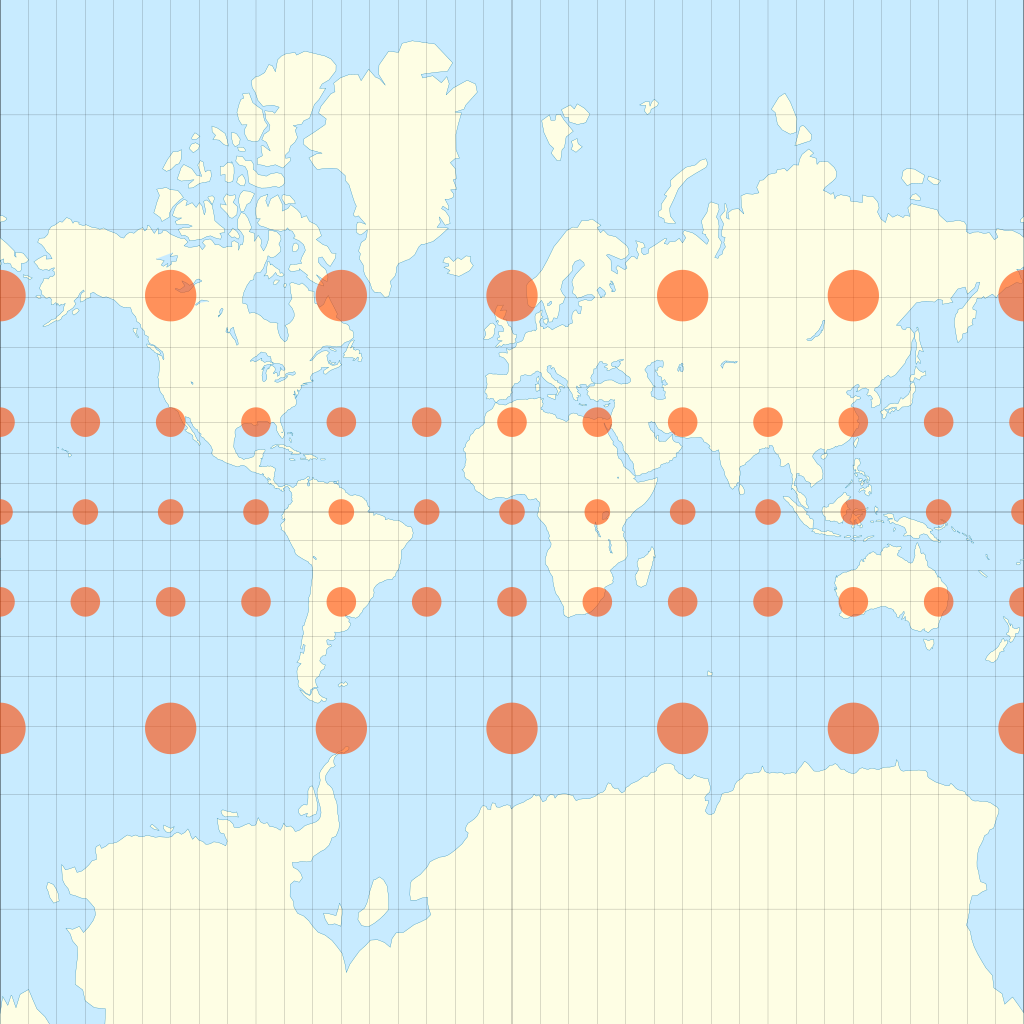}
\end{center}
 The Mercator projection has some excellent properties  preserving local directions and shapes, yet it also inflates the size of objects far from the equator. This is illustrated in the sample map (on the right) by red circles which in reality have the same area, yet on the map they have different sizes. Gauss demonstrated that every projection, no matter how clever, will always result in some sort of deformation like this.

As already mentioned, the name of Gauss is often linked to another
important development in geometry. In order to describe it, some preliminary
remarks are needed.
In the \emph{Elements}  Euclid gave five postulates for plane geometry. Postulates or axioms are statements which are believed to be obvious and require no proof.
 The first four of them assume that it is always possible
(1) to draw a straight line from any point to any point,
(2) to extend a finite linear segment continuously in a straight line
(3) to draw a circle with any centre and radius and that (4)  all right angles are equal to one another. The fifth postulate, called the parallel postulate  was far more complicated and not immediately obvious, namely 
\begin{quote}
If a straight line falling on two straight lines make the interior angles on the same side less than two right angles, the two straight lines, if produced indefinitely, meet on that side on which the angles are less than two right angles.
\end{quote}
Today we usually state the fifth postulate in the equivalent form given by
John Playfair (1748--1819): \emph{In a plane, given a line and a point not on it, at most one line parallel to the given line can be drawn through the point.}
 
Because the fifth postulate is not self-evident, for over two thousand years mathematicians tried to prove that postulate using the first four.
In the 19-th century it finally turned out that such efforts were futile because  the parallel postulate is independent of the remaining four. One can, therefore,
construct a self-consistent versions of geometry either with or without the fifth postulate, meaning that there is no \emph{one geometry}, but \emph{many geometries}. This remarkable and truly ground breaking result
has been discovered,  as it often happens,  independently and nearly simultaneously by two people,  the Hungarian mathematician   János Bolyai (1802--1860) and the Russian mathematician Nikolai Ivanovich Lobachevsky (1792--1856). Both Russia and Hungary
commemorated their fellow countrymen  with coins.

The Russian 1 Rouble copper-nickel coin commemorates the 200th anniversary of Lobachevsky's birthday. It is one of the first commmemorative coins minted in 1992 by the Russian Federation after the dissolution of the Soviet Union, bearing on the obverse  the new inscription ``Bank of Russia'' in place of the former ``CCCP''.
The reverse features a half-body portrait of Lobaczewski, 
\begin{center}
\includegraphics[width=5cm]{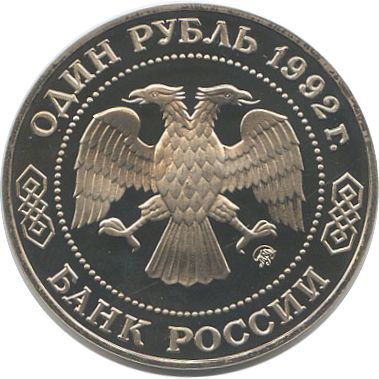}
\includegraphics[width=5cm]{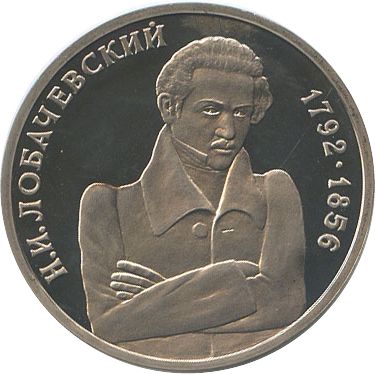}
\end{center}

The 200th anniversary of the birth of János Bolyai happened 10 years later, in 2002,
and  Hungary followed  suit by commemorating this event with a 3000 Forint silver coin.
Bolyai published his work on non-Euclidean geometry in 1831 as an appendix of a book
authored by his father. The reverse carries the full tile of this work (in Latin):
\begin{quote}
Appendix, Scientiam Spatii absolute veram exhibens; a veritate aut falsitate Axiomatis XI. Euclidei (a priori haud unquam decidenda) independentem; adjecta ad casum falsitatis quadratura circuli geometrica. 
\end{quote}
In English this translates as  ``Appendix, The absolute true Science of Space exhibited; independently of the XIth Euclidean Axiom (that can never be decided a priori) being true or false; for the case of being false the geometric quadrature of the circle is supplemented''.
The obverse shows an exact copy of one of the figures from the Appendix, namely
Figure 10, here reproduced on the left.
\begin{center}
\includegraphics[width=5cm]{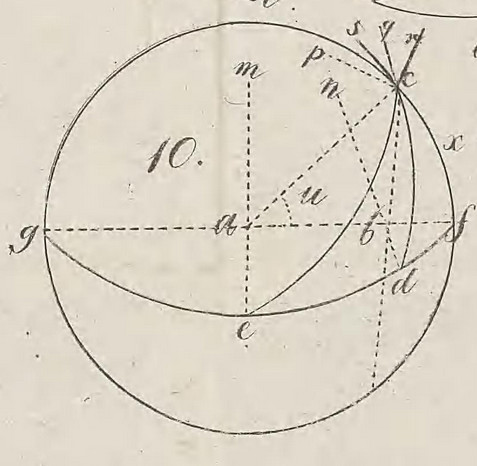}
\includegraphics[width=5cm]{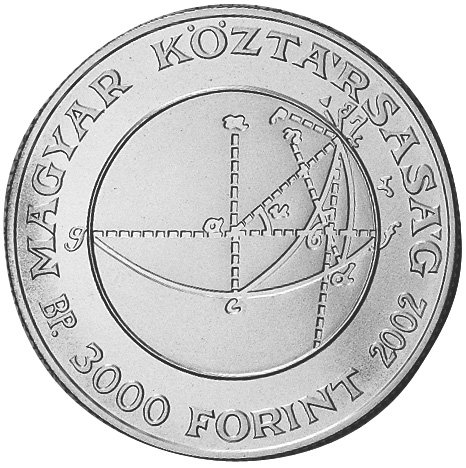}
\includegraphics[width=5cm]{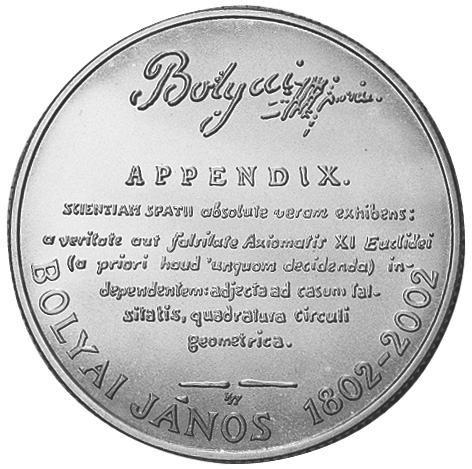}
\end{center}
During his lifetime Bolyai had not received the recognition he deserved. A copy of the Appendix work was sent to Gauss by Farkas Bolyai, father of
János. Gauss wrote in response that he was unable to praise János's work because it would be as if he praised himself since he had discovered the main results some years before. This was a devastating blow to Bolyai, even though Gauss had no claim to priority since he had never published the alleged findings.
Furthermore, after the death of Gauss no written proof of the aforementioned statement was found in any of his papers. 

The discoveries of Bolyai and Lobachevski fell into obscurity after their publication and were almost entirely forgotten until the 1860s, when the work of Bernhard Riemann on differential geometry
sparked a new wave of interest in non-Euclidean geometry. 
Riemann proposed a very general recipe for creating  geometry in
arbitrary space. His recipe requires only  the ability to measure distances between points, so
that given any two points you can talk about the shortest path between them, and
the ability to measure angles between intersecting curves that cross. 
The notion of distance allows to compute a quantity known as curvature,
which in the space with Euclidean geometry is exactly zero everywhere.
Although Riemann did not explicitly name them as such, the geometries of
Bolyai and Lobachevski are special cases of Riemannian geometry corresponding to constant negative curvature.
Today  such geometries are known as  Bolyai–Lobachevski 
or hyperbolic geometries.

After the second half of  the 19-th century geometry experienced
further intense growth. In addition to differential geometry, many new branches appeared, including  algebraic, finite and computational geometry. Some of these fields are very new and
we will have to wait a bit until we see the commemoration of their founders on coins.
The next opportunity, however, is very close.
In 2026 we will celebrate the 200th anniversary of the birth of Bernhard Riemann.
One can hope that his contributions to mathematics, including geometry, are also
honored with a  beautiful coin.

\providecommand{\href}[2]{#2}\begingroup\raggedright\endgroup

\begin{thebibliography}{1}

\bibitem{paper47}
H.~Fuk\'s, ``Mathematics on coins, part {I}: The tale of two queens and two
  towering figures,'' {\em The Canadian Numismatic Journal} {\bf 57} (2012),
  no.~5, 304--315, \href{http://xxx.lanl.gov/abs/arXiv:1601.07394}{{\tt
  arXiv:1601.07394}}.

\bibitem{paper58a}
H.~Fuk\'s, ``Mathematical formulae on coins. part 1: Einstein,'' {\em The
  Canadian Numismatic Journal} {\bf 61} (2016), no.~4, 167--169,
  \href{http://xxx.lanl.gov/abs/arXiv:1609.09485}{{\tt arXiv:1609.09485}}.

\bibitem{paper58b}
H.~Fuk\'s, ``Mathematical formulae on coins. part 2: from {P}ythagoras to the
  20th century,'' {\em The Canadian Numismatic Journal} {\bf 61} (2016), no.~5,
  208--213, \href{http://xxx.lanl.gov/abs/arXiv:1609.09485}{{\tt
  arXiv:1609.09485}}.

\end{thebibliography}
\end{document}